\documentclass[12pt]{amsart}
\usepackage{latexsym,enumerate}
\usepackage{amsmath,amsthm,amsopn,amstext,amscd,amsfonts,amssymb}
\usepackage{color}
\numberwithin{equation}{section}

\newtheorem{theorem}{Theorem}
\newtheorem{lemma}{Lemma}
\newtheorem{corollary}{Corollary}
\newtheorem{proposition}{Proposition}
\newtheorem{remark}{Remark}

\numberwithin{theorem}{section}
\numberwithin{corollary}{section}
\numberwithin{lemma}{section}
\numberwithin{definition}{section}
\numberwithin{proposition}{section}
\numberwithin{remark}{section}

\newcommand{\RR}{\mathbb R^N}

\newcommand{\medint}{-\kern  -,375cm\int}
\newcommand{\dint}{\displaystyle\int}
\newcommand{\dsum}{\displaystyle\sum}

\begin{document}

\title[An isoperimetric inequality for Gauss--like product measures]
{An isoperimetric inequality \\for Gauss--like product measures}

\author{ F. Brock$^1$ - F. Chiacchio$^2$ - A. Mercaldo$^2$}
\thanks{}
\date{}

\begin{abstract}
\noindent This paper deals with various questions related to the isoperimetic problem
 for smooth positive measure $d\mu = \varphi(x)dx$, with $x \in \Omega \subset \RR$.
Firstly we find some necessary conditions on the density of the measure $ \varphi(x)$  
that render the intersection of half spaces with $\Omega$ a minimum in the isoperimetric problem. 
We then identify the unique isoperimetric set for a wide class of factorized finite measures.
These results are finally used in order to get sharp  
inequalities in weighted Sobolev spaces and a comparison result for solutions
to boundary value problems for degenerate elliptic equations.
 \bigskip


\end{abstract}

\maketitle

\setcounter{footnote}{1} 
\footnotetext{Leipzig University, Department of Mathematics, Augustusplatz, 04109 Leipzig,
Germany, e-mail: brock@math.uni-leipzig.de}

\setcounter{footnote}{2} 
\footnotetext{Dipartimento di Matematica e Applicazioni \textquotedblleft R.
Caccioppoli\textquotedblright , Universit\`{a} degli Studi di Napoli
``Fe\-derico II", Complesso Monte S. Angelo, via Cintia, 80126 Napoli,
Italy, e-mails: francesco.chiacchio@unina.it, mercaldo@unina.it}

\section{Introduction}
This paper deals with relative isoperimetric inequalities in the setting of 
\textsl{manifold with density}. 
 More precisely let $\Omega $ be a Lebesgue measurable set in $\mathbb{R}^{N}$  and let 
$\mu $ be a positive finite measure on $\Omega $ given by 
\begin{equation}
d\mu (x) =\varphi (x)\,dx,  \label{dmu}
\end{equation}
where $\varphi$ is a positive function in 
$ C^{0}(\Omega )$ and $\mu (\Omega )<+\infty $.
For any Borel measurable subset $M $ of $\Omega $, the $\mu $--perimeter of $M$
relative to $\Omega $ is given by 
\begin{equation}
P_{\mu }(M,\Omega ):=\sup \left\{
\int_{M}\mbox{div}\,\left( \varphi  v \right)
\,dx:\,  v \in C_{0}^{1}(\Omega ,\mathbb{R}^{N}),\,| v |\leq 1\mbox{ in }\,\Omega
\right\} .
\end{equation}

\noindent As well known, the above \textsl{distributional definition} of weighted perimeter is equivalent to the following 

\begin{equation}
P_{\mu }(M,\Omega )=
\left\{ 
\begin{array}{ccc}
\dint_{\partial M \cap \Omega}\varphi (x)\,{\mathcal{H}}
_{N-1}(dx) & 
\mbox{ if } &   \partial M \cap  \Omega  \mbox{ is } (N-1)-\mbox{rectifiable } \\ 
&  &  \\ 
+ \infty \qquad  & \mbox{ otherwise. } & 
\end{array}
\right. 
\end{equation}

\noindent The function $\varphi (x)$, which appears both in the volume and in the perimeter, is called the \textsl{density}.
We say that a set is isoperimetric 
or solves the isoperimetric problem relative to 
$\Omega$ if it minimizes the weighted perimeter $P_{\mu }(M,\Omega )$ among
all the sets $M \subset \Omega $
with fixed weighted  measure $\mu (M)$.
This subject has attracted a growing interest starting from the papers by Sudakov-Tsirel'son and Borell 
(see \cite{Su} and \cite{Bo}) 
on the isoperimetric problem for Gaussian density, where it turned out that the isoperimetric  set is a half-space. 
Since then the isopermetric problem has been solved for various class of weights
(see, e.g.,  \cite{E}, \cite{Bobkov}, 
\cite{MS}, \cite{CK}, \cite{BCM2008}, \cite{BBMP}, \cite{Rosales},  \cite{BCM2012}, 
\cite{RCBM},  \cite{BMP}, \cite{BBC}, \cite{CR-O},
\cite{BCM2014}, \cite{CFMP},  \cite{FPR}  and \cite{Chambers}). 
Clearly such a bibliography if far from being exhaustive.

We are interested in two types of questions. Firstly we find some necessary conditions on smooth positive measures  
$\mu$ that render the intersection of half spaces with $\Omega $ a minimum in the relative isoperimetric problem.
Among other things, we show that the weight function $\varphi$ must be in separated form 
$$
\varphi (x)=\rho (x_1,\ldots ,x_{N-1}) \sigma (x_{N})
$$
for some positive functions $\rho$ and $\sigma$ (see Theorem \ref{Theorem1} in the next section).
In Theorem \ref{Theorem2}, our main result, we identify the unique isoperimetric set for a wide class of factorized finite measures.

\noindent  In order to state this last  result we need some notation.
For $i=1,\ldots ,N-1,$ ($N\geq 2$), let $-\infty \leq a_{i}<b_{i}\leq +\infty $, and let $A_{i}\in C^{1}(a_{i},b_{i})$ be  
real functions  such that
\begin{equation}
A_{i}^{\prime }(x)\geq 1\quad \mbox{ on }\ (a_{i},b_{i}),  \label{keyineq}
\end{equation}
and
\begin{equation*}
\lim_{x\rightarrow a_{i}^{+}}A_{i}(x)=-\infty 
\quad  \mbox{and}  \quad
\lim_{x\rightarrow b_{i}^{-}}A_{i}(x)=+\infty .
\end{equation*}
Further,  let
 \begin{equation*}
S^{\prime }:=(a_{1},b_{1})\times \cdots \times
(a_{N-1},b_{N-1}) \quad  \mbox{and}  \quad
S:=S^{\prime }\times \mathbb{R}, 
 \end{equation*}
and, finally,  let $\mu $ be the measure on $S$, given by 
\begin{equation}
d\mu (x) :=\varphi (x)\,dx^{\prime }\,dx_{N} = \varphi (x)\,dx,
\label{defmu}
\end{equation}
where 
\begin{equation}
\varphi (x):=\mbox{exp}\,\left\{ -\sum_{i=1}^{N-1}\frac{A_{i}(x_{i})^{2}}{2}-
\frac{x_{N}^{2}}{2}\right\} \prod_{i=1}^{N-1}A_{i}^{\prime }(x_{i}),\ \
x\in S.  \label{defphi}
\end{equation}

\noindent  If $\lambda \in \mathbb{R}$, let $S_{\lambda }$ be the intersection of the halfspace 
$\{x_{N}>\lambda \}$ with $S$, that is, 
\begin{equation*}
S_{\lambda }= 
S^{\prime} \times ( \lambda ,\infty ).
\end{equation*}

\begin{theorem}
\label{Theorem2}
Let $M$ be a Lebesgue measurable subset of $S$ and fix $\lambda $ such that
\begin{equation}
\mu (M)=\mu (S_{\lambda }).  \label{measM}
\end{equation}
Then 
\begin{equation}
P_{\mu }(M,S)\geq P_{\mu }(S_{\lambda },S).  \label{Pmuisop}
\end{equation}
Moreover equality holds in \eqref{Pmuisop} if and only $M=S_{\lambda }$.
\end{theorem}

As we will show in Corollary \ref{Corollary1}, the conclusion of Theorem \ref{Theorem2} holds for measures $\mu $ of the
type 
\begin{equation}
d\mu (x):=\exp \left\{ -\frac{|x|^{2}}{2}-\sum_{i=1}^{N-1}B_{i}(x_{i})\right
\} \,dx,  \label{measure2}
\end{equation}
where $B_{i}\in C^{2}(a_{i},b_{i})$ with $B_{i}^{\prime \prime }(x_{i})\geq 0
$ on $(a_{i},b_{i})$, ($i=1,\ldots ,N-1$). 
\\[0.2cm]
Our isoperimetric inequality Theorem \ref{Theorem2} is proved in Section 4. It generalizes the results contained in \cite{Rosales} in two directions. 
We consider more general factorized
perturbations of the Gaussian measure and we allow these perturbations to affect not just one but $N-1$ variables. 
Note that in view of Lemma \ref{Lemma1} and the remark following Corollary \ref{Corollary1} in Section 4, the weight function in  (\ref{defphi}) is indeed 
more general than the one of (\ref{measure2}).
The main ingredient  in the proof of Theorem \ref{Theorem2} consists in using a map that coincides with the optimal 
 transport Brenier map and that pushes 
the measure $d\mu$ forward to the Gaussian measure. We explicitly remark that the relevant property of the gradient 
of such a  map is proved by means of elementary and self-contained tools.
While in \cite{Rosales} the analogous question is faced by using a result by Cafarelli (see \cite{Ca}).
 An important example for (\ref{measure2}) is given by $a_{i}=0$, $b_{i}=+\infty $, $
B_{i}(x_{i})=-k_{i}\log x_{i}$ with $k_{i}\geq 0$, ($i=1,\ldots ,N-1$), that
is, 
\begin{equation}
d\mu (x)=\exp \left\{ -\frac{|x|^{2}}{2}\right\}
\prod_{i=1}^{N-1}x_{i}^{k_{i}}\,dx.  \label{specmeas}
\end{equation}
Finally, in  Section 5, using a kind of symmetrization, related to the isoperimetric inequalities that we have proved,
we give some sharp apriori  bounds to the solutions of a class of elliptic second order Pde's (see Theorem \ref{Theorem3}).

\section{Necessary conditions}

Below we will introduce some weighted spaces: 
for $p\in [1, +\infty ]$, let $L^p (\Omega , d\mu )$ 
be the standard weighted $L^p$-space (corresponding to the weight $d\mu =\varphi (x)dx $). 
By $W^{1,2} (\Omega , d\mu ) $ we denote the weighted Sobolev space,
$$
W^{1,2} (\Omega , d\mu ) := \{ u \in W^{1,1} _{loc} (\Omega ):\, u,|\nabla u| \in L^2 (\Omega , d\mu \} .
$$

 We begin our analysis with some necessary conditions on smooth positive
finite measures $\mu $ that render the intersection of half spaces with 
$\Omega $ a minimum in the relative isoperimetric problem. 
\newline
Following \cite{RCBM} and \cite{Rosales}, we introduce the notion of stationarity and stability of sets. 
Let $\Omega $ be a smooth set with boundary $\Sigma $ and inward unit normal vector $\bold\nu$. 
We consider a one-parameter variation 
$\{ \phi _t \} _{|t|<\varepsilon }: \mathbb{R} ^N \to \mathbb{R}^N $ 
with associated infinitesimal vector field $X= d\phi _t / dt $ with  normal component $u = \langle X,\bold \nu\rangle $. 
Let $\Omega _t = \phi _t (\Omega )$ 
and $\Sigma _t = \phi _t (\Sigma )$. The volume and perimeter functions of the variation are $V(t):= \mu (\Omega _t) $ 
and $P(t):= P_{\mu } (\Omega _t )$, respectively. 
We say that a given variation $\{ \phi _t \} _t $ 
{\sl preserves volume} if $V(t)$ is constant for any small $|t|$.
We say that $\Omega $ is {\sl stationary} if $P'(0)=0$ for any volume-preserving variation. 
Obviously any isoperimetric region is also stationary. Finally, we say that $\Omega $ is {\sl stable} 
if it is stationary and if $P''(0) \geq 0$ for any volume-preserving variation of $\Omega $. 
We note that the first and second variation 
of the volume and perimeter, $V'(0)$, $P'(0)$, $V''(0)$ and $P''(0)$, respectively, were given in \cite{RCBM}. 
\newline
The following notation for points and the gradient in $\mathbb{R}^{N}$ 
will be in force throughout the paper
$$x=(x^{\prime },x_{N}), \hspace{0.2cm}  x^{\prime }=(x_{1},\ldots ,x_{N-1}) \in
\mathbb{R}^{N-1}, \hspace{0.2cm} x_{N}\in \mathbb{R}
$$
and
$$
\nabla =(\nabla ^{\prime },\partial /\partial
x_{N}),  \hspace{0.2cm}  \nabla ^{\prime }=(\partial /\partial x_{1},\ldots ,\partial
/\partial x_{N-1}).
$$
If $\Omega ' $ is a domain in 
$\mathbb{R}^{N-1}$ we set 
\begin{equation}
\Omega _{\lambda }:=\{(x^{\prime },x_{N}):\,x_{N}>\lambda ,\,x^{\prime }\in \Omega '
\},\quad \lambda \in \mathbb{R}.  \label{slambda}
\end{equation}
Our first result is

\begin{theorem}
\label{Theorem1}
 Let $\Omega :=\Omega ' \times \mathbb{R}$ where $\Omega ' $ 
is a domain in 
$\mathbb{R}^{N-1}$ with Lipschitz boundary, and let $\mu $ be a measure on $\Omega $
given by 
\begin{equation}
d\mu (x) =\varphi (x)\,dx,\quad x\in \Omega ,  
\label{dmu1}
\end{equation}
where $\varphi \in C^{1}(\Omega )$ and $\varphi (x)>0$ on $\Omega $. \newline
(i) If $\Omega _{\lambda }$ is stationary in the relative isoperimetric problem for 
$\mu $ and $\Omega $, for every $\lambda \in \mathbb{R}$, then 
\begin{equation}
\varphi (x)=\rho (x^{\prime })\sigma (x_{N})
\quad 
\forall x\in \Omega ,
\label{decomp}
\end{equation}
where $\rho \in C^{1}(\Omega ' )$ and $\sigma \in C^{1}(\mathbb{R})$ are positive. 
\newline
(ii) If $S_{\lambda }$ is stable in the relative isoperimetric problem for $
\mu $ and $\Omega $, for every $\lambda \in \mathbb{R}$, then 
\begin{equation}
\kappa _{1}\geq \tau ,  \label{stabineq}
\end{equation}
where 
\begin{equation}
\tau :=\sup \left\{ \left( 
\frac{\sigma ^{\prime }(t)}{\sigma (t)}
\right)
^{2}-\frac{\sigma ^{\prime \prime }(t)}{\sigma (t)}:\,
t\in \mathbb{R}
\right\} ,  \label{tau}
\end{equation}
and 
\begin{equation}
\kappa _{1}:=\inf \left\{ \frac{\int_{\Omega ' }|\nabla ^{\prime }v | ^2 \rho
\,dx^{\prime }}{\int_{\Omega ' }v^{2}\rho \,dx^{\prime }}:\,
v\in
W^{1,2}(\Omega ' ,\rho dx' ),\,
v\not\equiv 0\right\} .
\end{equation}

\end{theorem}

\begin{remark}\rm
{\bf (a)} Observe that $\kappa _{1}$ is the first nontrivial eigenvalue of the
Neumann problem 
\begin{equation}
\left\{ 
\begin{array}{ccc}
-\dsum_{k=1} ^{N-1} \frac{\partial }{\partial x_k}  \left( \rho  \frac{\partial}{\partial x_k }  u\right) =
\kappa \rho u & 
\mbox{ in } & \Omega '  \\ 
&  &  \\ 
\displaystyle{\frac{\partial u}{\partial \mathbf{n}} }=0\qquad  & \mbox{ on } & \partial
\Omega ' ,
\end{array}
\right. 
\end{equation}
where $u\in W^{1,2}(\Omega ',\rho dx' )$, 
and $\mathbf{n}$ is the exterior unit normal to $\partial \Omega ' $. \newline
{\bf (b)} Let $\Omega =\mathbb{R}^{N}$. If $\sigma (t)=e^{-ct^{2}}$ for some $c>0$
, then $\tau =2c$, and if $\rho (x^{\prime })=e^{-c\left\vert x^{\prime
}\right\vert ^{2}}$, then also $\kappa _{1}=2c$, (see \cite{Tay}, p.105 ff.), 
so that condition (\ref{stabineq}) is satisfied for Gauss
measures, $\varphi (x)=e^{-c|x|^{2}}$, ($x\in \mathbb{R}^{N}$). 
\end{remark}

{\sl Proof of Theorem \ref{Theorem1}:}  Proceeding similarly as in \cite{BCM2008}, we define
volume--preserving perturbations from $\Omega _{\lambda }$. Let $u\in C^{2}(\Omega '
)$. Then the Implicit Function Theorem tells us that there exists a number 
$\varepsilon _{0}>0$ and a function $s\in C^{2}(-\varepsilon
_{0},\varepsilon _{0})$ with $s(0)=0$, such that 
\begin{equation*}
\Omega _{\lambda }(\varepsilon ):= 
\{(x^{\prime },x_{N}):\,x_{N}>u(x^{\prime},\varepsilon ), x^{\prime }\in \Omega \}, 
\end{equation*}
where
\begin{equation*}
 u(x^{\prime },\varepsilon ) := \lambda + \varepsilon
u(x^{\prime })+s(\varepsilon ),
\end{equation*}
and $\Omega _{\lambda }$ have the same $\mu $--measure, that is 
\begin{equation}
\mu (\Omega _{\lambda }(\varepsilon ))=\int_{\Omega ' }\int_{u(x^{\prime
},\varepsilon )}^{+\infty }\varphi (x^{\prime },t)\,dt\,dx^{\prime }=\mu
(\Omega _{\lambda }).  
\label{measpreserve1}
\end{equation}
This implies 
\begin{equation}
0=\frac{\partial }{\partial \varepsilon }\mu (\Omega _{\lambda }(\varepsilon
))=-\int_{\Omega ' }\varphi (x^{\prime },u(x^{\prime },\varepsilon ))\left(
u(x^{\prime })+s^{\prime }(\varepsilon )\right) \,dx^{\prime },
\label{measpreserve2}
\end{equation}
for $|\varepsilon |<\varepsilon _{0}$. Writing $s_{1}:=s^{\prime }(0)$ and $
s_{2}:=s^{\prime \prime }(0)$, we obtain 
\begin{equation}
0=-\frac{\partial }{\partial \varepsilon }\mu ((\Omega _{\lambda }(\varepsilon ))
\Big|_{\varepsilon =0}=\int_{\Omega ' }\varphi (x^{\prime }, \lambda )(u(x^{\prime
})+s_{1})\,dx^{\prime }.  \label{measpreserve3}
\end{equation}
Further, we have 
\begin{equation}
P_{\mu }(\Omega _{\lambda }(\varepsilon ),\Omega )=\int_{\Omega ' }
\varphi
(x^{\prime },u(x^{\prime },\varepsilon ))
\sqrt{1+\varepsilon ^{2}
|\nabla ^{\prime } u(x^{\prime } )|^2 }
\, dx^{\prime },  \label{perim}
\end{equation}
so that 
\begin{eqnarray}
&&
\frac{\partial }{\partial \varepsilon }P_{\mu }(\Omega _{\lambda }(\varepsilon
),\Omega )  \notag  \label{derPmu1} \\
&=&
\int_{\Omega ' }\left\{ \varphi _{x_{N}}(x^{\prime },u(x^{\prime
},\varepsilon ))
\left( u(x^{\prime })+s^{\prime }(\varepsilon )\right) 
\sqrt{
1+\varepsilon ^{2}|\nabla ^{\prime }u(x^{\prime} )| ^2 } +
\right.   \notag \\
&& \left. +\varepsilon \varphi (x^{\prime },u(x^{\prime },\varepsilon
))
\left( 1+\varepsilon ^2 |\nabla ^{\prime }u(x^{\prime} )|^2 \right)
^{-1/2}
|\nabla ^{\prime }u(x^{\prime})|^2
\right\} \, dx^{\prime }.
\end{eqnarray}
Now assume that $\Omega_{\lambda }$ is stationary for every 
$\lambda \in \mathbb{R}$. Then (\ref{derPmu1}) gives 
\begin{equation}
\label{derPmu2}
0=
\frac{\partial }{\partial \varepsilon }P_{\mu }(\Omega _{\lambda}(\varepsilon ),\Omega )\Big|_{\varepsilon =0}
\hspace{.1 cm}
=\int_{\Omega ' }\varphi _{x_{N}}(x^{\prime },\lambda )(u(x^{\prime
})+s_{1})\,dx^{\prime }.
\end{equation}
This, together with (\ref{measpreserve3}) implies that $\int_{\Omega ' }\varphi
_{x_{N}}(x^{\prime },\lambda )v(x^{\prime })\,dx^{\prime }=0$ for all
functions $v\in C^{1}(\Omega ' )$ satisfying $\int_{\Omega ' }\varphi (x^{\prime
},\lambda )v(x^{\prime })\,dx^{\prime }=0$. Then the Fundamental Lemma in the Calculus
of Variations tells us that there is a number $k=k(\lambda )\in \mathbb{R} $ such that 
\begin{equation}
\varphi _{x_{N}}(x^{\prime },\lambda )=k(\lambda )\varphi (x^{\prime },\lambda
)\quad \forall x^{\prime }\in \Omega ',  \label{klambda}
\end{equation}
which implies (\ref{decomp}). Hence we have from (\ref{measpreserve2}) 
\begin{equation}
0=\int_{\Omega '}\rho (x^{\prime })\sigma (u(x^{\prime },\varepsilon ))\left(
u(x^{\prime })+s^{\prime }(\varepsilon )\right) \,dx^{\prime }. \label{m2}
\end{equation}
For $\varepsilon =0$ this yields 
\begin{equation}
\label{m3}
0 =\int_{\Omega '}\rho (x^{\prime }) \left(
u(x^{\prime })+s_1 \right) \, dx^{\prime }. 
\end{equation}
Differentiating (\ref{m3}) we obtain for $\varepsilon =0$,
\begin{equation}
0=\int_{\Omega ' }\rho (x^{\prime })\left[ \left( u(x^{\prime })+s_{1}\right)
^{2}\sigma ^{\prime }(\lambda )+s_{2}\,\sigma (\lambda )\right] \,dx^{\prime
}.  \label{measpreserve4}
\end{equation}
(ii) Next assume that $\varphi \in C^{2}(\Omega )$, and that $\Omega _{\lambda }$ is
stable for every $\lambda \in \mathbb{R}$. Then $\rho \in C^{2}(\Omega ' )$ and $\sigma
\in C^{2}(\mathbb{R})$. First, by (\ref{decomp}) and (\ref{derPmu1}) we have 
\begin{eqnarray*}
\frac{\partial }{\partial \varepsilon }P_{\mu }(\Omega _{\lambda }(\varepsilon
), \Omega ) &=&
\int_{\omega }\rho (x^{\prime })
\left\{ 
\sigma ^{\prime
}(u(x^{\prime },\varepsilon ))
\left( 
u(x^{\prime })+s^{\prime }(\varepsilon
)\right) \sqrt{1+\varepsilon ^{2}|\nabla ^{\prime }u(x^{\prime })|^2 }+\right. 
\\
&&\left. +\varepsilon \sigma (u(x^{\prime },\varepsilon ))\left(
1+\varepsilon ^{2}|\nabla ^{\prime }u(x^{\prime })|^2 \right) ^{-1/2}
|\nabla
^{\prime }u(x^{\prime })|^2 \right\} \, dx^{\prime }.
\end{eqnarray*}
Differentiating this gives 
\begin{eqnarray*}
0 &\leq &\frac{\partial ^{2}}{\partial \varepsilon ^{2}}P_{\mu }(\Omega _{\lambda
}(\varepsilon ))\Big|_{\varepsilon =0} \\
&=&\int_{\Omega ' }\rho (x^{\prime })\left\{ \sigma ^{\prime \prime }(\lambda
)\left( u(x^{\prime })+s_{1}\right) ^{2}+\sigma ^{\prime }(\lambda
)s_{2}+\sigma (\lambda )|\nabla ^{\prime }u(x^{\prime })|^2 \right\}
\,dx^{\prime }.
\end{eqnarray*}
In view of (\ref{measpreserve4}) we obtain 
\begin{equation*}
\int_{\Omega '}\rho (x^{\prime })|\nabla ^{\prime }u(x^{\prime }) |^2 \, dx^{\prime
}\geq \left( \left( \frac{\sigma ^{\prime }(\lambda )}{\sigma (\lambda )}
\right) ^{2}-\frac{\sigma ^{\prime \prime }(\lambda )}{\sigma (\lambda )}
\right) \int_{\Omega ' }\rho (x^{\prime })\left( u(x^{\prime })+s_{1}\right)
^{2}\,dx^{\prime }.
\end{equation*}
Hence (\ref{stabineq}) follows by (\ref{m3}). $\hfill \Box $

\section{A onedimensional auxiliary result}

Let $I:=(a,b)$, where $-\infty \leq a<b\leq +\infty $ and $B\in C^{2}(I)$
with $B^{\prime \prime }(x)\geq 0$ on $I$. Further, let 
\begin{eqnarray}
\nonumber
E(y):= &&\frac{1}{\sqrt{2\pi }}\int_{-\infty }^{y}e^{-t^{2}/2}\,dt,\quad
y\in \mathbb{R}, 
\\
\label{defc} 
c:= &&\frac{1}{\int_{a}^{b}e^{-t^{2}/2-B(t)}\,dt},
\end{eqnarray}
and let $A\in C^{3}(I)$ be given by 
\begin{equation}
A(x):=E^{-1}\left( c\int_{a}^{x}e^{-t^{2}/2-B(t)}\,dt\right) ,\quad x\in I.
\label{defA}
\end{equation}
Note that the convexity of $B$ ensures the convergence of the integrals on
the right-hand sides of (\ref{defc}) and (\ref{defA}), and that 
\begin{equation}
\label{formulaA}
e^{-A(x)^{2}/2}A^{\prime }(x)=c\sqrt{2\pi }e^{-x^{2}/2-B(x)},
\quad x \in I, 
\end{equation}
\begin{equation}
\label{Alimits}
\lim_{x\rightarrow a^{+}}A(x) =-\infty ,
\quad \lim_{x \rightarrow b^{-}}A(x) = +\infty .
\end{equation}
We also emphasize that the map $A$ of (\ref{defA}) 
coincides with the optimal transport Brenier map pushing the measure 
$$
d\mu_1 (x) := c e^{-x^2 /2 -B(x)} \, dx ,  
$$
defined on $(a,b)$, 
forward to the one-dimensional Gauss measure,
$$
d\gamma _1 (y) := \frac{1}{\sqrt{2\pi } } e^ {-y^2 /2 } \, dy, 
$$
(see \cite{Ca}, Thm. 1 and 2). 
Hence we can use a result of \cite{Ca}, Thm.11, to obtain the following Lemma \ref{Lemma1}. 
For the convencience of the reader, we include an elementary proof.

{
\begin{lemma}\label{Lemma1}
 Under the assumptions above, 
\begin{equation}
A^{\prime }(x)\geq 1,\qquad  x\in I.  \label{A'geq1}
\end{equation}
\end{lemma}

{\sl Proof:} We first claim: 
\begin{equation}
\label{localmin}
\mbox{If $A'$ has a local minimum at $x_0 \in I$, 
 then $A'(x_0 )\geq 1 $.}
\end{equation}
Assume that $A^{\prime }$ has a local minimum at $x_{0}\in I$ and $A^{\prime
}(x_{0})<1$. Identity (\ref{formulaA}) gives 
\begin{eqnarray*}
A^{\prime \prime }(x) &=&A^{\prime }(x)\left( A(x)A^{\prime }(x)-x-B^{\prime
}(x)\right) \quad \mbox{and} \\
A^{\prime \prime \prime }(x) &=&A^{\prime \prime }(x)\left( A(x)A^{\prime
}(x)-x-B^{\prime }(x)\right) + \\
&&+A^{\prime }(x)\left( [A^{\prime 2}+A(x)A^{\prime \prime }(x)-1-B^{\prime
\prime }(x)\right)
\end{eqnarray*}
on $I$. Since $A^{\prime \prime }(x_{0})=0\leq A^{\prime \prime \prime
}(x_{0})$, this implies 
\begin{eqnarray*}
0 &\leq &A^{\prime }(x_{0})\left( [A^{\prime }(x_{0})]^{2}-1-B^{\prime
\prime }(x_{0})\right) \\
&\leq &A^{\prime }(x_{0})\left( [A^{\prime }(x_{0})]^{2}-1]\right) <0,
\end{eqnarray*}
a contradiction. Hence (\ref{localmin}) holds. \newline
This implies that \eqref{A'geq1} holds for points inside $I$. 
It remains to show that 
\begin{equation}
\label{lima}
\liminf_{x\rightarrow a^{+}}A^{\prime }(x) \geq 1
\end{equation}
and   
\begin{equation}
\label{limb}
\liminf_{x\rightarrow b^{-}}A^{\prime }(x) \geq 1.
\end{equation}
We only show (\ref{lima}). 
The proof of (\ref{limb}) is similar and is left
to the reader. \newline
Assume by absurd that 
$$
\liminf_{x\rightarrow a^{+}}A^{\prime }(x)=:L<1.
$$
 By (\ref{localmin}) this implies that $\lim_{x\rightarrow a^{+}}A^{\prime }(x)$ exists
and 
\begin{equation}
\lim_{x\rightarrow a^{+}}A^{\prime }(x)=L.  \label{limA'}
\end{equation}
By (\ref{formulaA}) and (\ref{Alimits}), this means in particular that 
$ \lim_{x\rightarrow a^{+}}(x^{2}/2+B(x))=+\infty $. 
In view of the convexity of
 $B$, we deduce 
\begin{equation}
\lim_{x\rightarrow a^{+}}(x+B^{\prime }(x))=-\infty .  \label{limB'}
\end{equation}
Then the generalized Mean Value Theorem tells us that for every $x\in (a,b)$
there exists a number $y\in (a,x)$ such that 
\begin{eqnarray}
A^{\prime }(x) &=&c\sqrt{2\pi }\frac{e^{-x^{2}/2-B(x)}}{e^{-A(x)^{2}/2}} 
\notag  \label{xy} \\
&=&c\sqrt{2\pi }\frac{[e^{-y^{2}/2-B(y)}]^{\prime }}{[e^{-A(y)^{2}/2}]^{
\prime }}  \notag \\
&=&c\sqrt{2\pi }\frac{(y+B^{\prime }(y)) e^{-y^{2}/2-B(y)}}{A(y)A^{\prime } (y) e^{ -A(y)^{2}/2}
}  \notag \\
&=&\frac{y+B^{\prime }(y)}{A(y)}.
\end{eqnarray}
In view of (\ref{Alimits}) and (\ref{limB'}) and since $L<1$, we find a
strictly decreasing sequence $\{x_{n}\}$ with $x_{n}\rightarrow a$ such that 
\begin{equation}
\frac{1-\frac{A(x_{n})}{A(x_{n+1})}}{1-\frac{x_{n}+B^{\prime }(x_{n})}{
x_{n+1}+B^{\prime }(x_{n+1})}}\geq \frac{L(L+1)}{2}  \label{ineq1}
\end{equation}
Using once more the generalized Mean Value Theorem and (\ref{ineq1}) we find
another sequence $\{y_{n}\}$ with $x_{n+1}<y_{n}<x_{n}$ and such that 
\begin{eqnarray*}
\frac{x_{n+1}+B^{\prime }(x_{n+1})}{A(x_{n+1})} &=&\frac{x_{n}+B^{\prime
}(x_{n})-x_{n+1}-B^{\prime }(x_{n+1})}{A(x_{n})-A(x_{n+1})}\cdot \frac{1-
\frac{A(x_{n})}{A(x_{n+1})}}{1-\frac{x_{n}+B^{\prime }(x_{n})}{
x_{n+1}+B^{\prime }(x_{n+1})}} \\
&=&\frac{1+B^{\prime \prime }(y_{n})}{A^{\prime }(y_{n})}\cdot \frac{1-\frac{
A(x_{n})}{A(x_{n+1})}}{1-\frac{x_{n}+B^{\prime }(x_{n})}{x_{n+1}+B^{\prime
}(x_{n+1})}} \\
&\geq &\frac{1+B^{\prime \prime }(y_{n})}{A^{\prime }(y_{n})}\cdot \frac{
L(L+1)}{2} \\
&\geq &\frac{L(L+1)}{2A^{\prime }(y_{n})}\longrightarrow \frac{L+1}{2},\quad 
\mbox{as $\ n\to \infty $.}
\end{eqnarray*}
Hence 
\begin{equation*}
\limsup_{x\rightarrow a^{+} }\frac{x+B^{\prime }(x)}{A(x)}\geq \frac{L+1}{2}.
\end{equation*}
By (\ref{xy}) this means that also 
\begin{equation*}
\limsup_{x\rightarrow a^{+}}A^{\prime }(x)\geq \frac{L+1}{2}>L,
\end{equation*}
contradicting (\ref{limA'}). It follows that $L\geq 1$. $\hfill \Box $
\\[0.1cm]

\section{The isoperimetric inequality}
In this section we prove our main result Theorem \ref{Theorem2}. Let $\gamma _{N}$ denote the $N$-dimensional Gauss measure on $\mathbb{R}^{N}
$, given by 
\begin{equation}
d\gamma _{N}(x):=(2\pi )^{-N/2}e^{-|x|^{2}/2}dx.  \label{gaussmeas}
\end{equation}
For any $U$ Lebesgue measurable subset of  $\mathbb{R}^{N}$, let $P_{\gamma _{N}}(U)$ denotes its 
Gaussian perimeter.

\vspace{0.3cm}

\noindent  \textsl{Proof of Theorem \ref{Theorem2}:}  Define a diffeomorphism $T$ between $S$ and $\mathbb{R}^{N}$, by 
$$T(x^{\prime },x_{N}):=(\tilde{T}(x^{\prime }),x_{N}),$$ 
where 
$$
\tilde{T} (x^{\prime })
:=(A_{1}(x_{1}),\ldots ,A_{N-1}(x_{N-1})), \hspace{0.2cm} (x^{\prime},x_{N})\in S,
$$ 
and let 
\begin{equation*}
H_{\lambda }:=\{(x^{\prime },x_{N}):\,x_{N}>\lambda ,\,x^{\prime }\in 
\mathbb{R}^{N-1}\}.
\end{equation*}
Clearly we have 
\begin{eqnarray*}
&&T(S_{\lambda })=H_{\lambda }\ \mbox{ and } \\
&&\mu (M)=(2\pi )^{N/2}\gamma _{N}\left( T(M)\right) ,\ \mu (S_{\lambda
})=(2\pi )^{N/2}\gamma _{N}\left( T(H_{\lambda })\right) .
\end{eqnarray*}
Hence (\ref{measM}) together with the isoperimetric inequality in Gauss
space yields 
\begin{equation}
P_{\gamma _{N}}\left( T(M)\right) \geq P_{\gamma _{N}}(H_{\lambda }).
\label{isopgauss}
\end{equation}
Since also 
\begin{eqnarray*}
P_{\mu }(S_{\lambda },S) &=&\int_{S^{\prime }}\varphi (x^{\prime },\lambda
)\,dx^{\prime } \\
&=&\int_{\mathbb{R}^{N-1}}\mbox{exp}\,\left\{ -\frac{\left\vert x^{\prime
}\right\vert ^{2}+\lambda ^{2}}{2}\right\} \,dx^{\prime } \\
&=&(2\pi )^{N/2}P_{\gamma _{N}}(H_{\lambda }),
\end{eqnarray*}
it remains to show that 
\begin{equation}
P_{\mu }(M,S)\geq (2\pi )^{N/2}P_{\gamma _{N}}\left( T(M)\right) .
\label{perineq}
\end{equation}
To prove (\ref{perineq}), we first consider the case that $\Sigma $ is an open subset of $S\cap \partial
M$ given in the form 
\begin{equation}
\Sigma =\{(x^{\prime },u(x^{\prime })):\,x^{\prime }\in \Sigma ^{\prime }\},
\label{defsigma}
\end{equation}%
where $u\in C^{1}(\Sigma ^{\prime })$ and $\Sigma ^{\prime }$ is an open
subset of $S^{\prime }$. We write $y^{\prime }\equiv \tilde{T}(x^{\prime })$
, $v(y^{\prime }):=u(x^{\prime })$, ($x^{\prime }\in \Sigma ^{\prime }$),
and $\tilde{T}(\Sigma ^{\prime }):=\{\tilde{T}(x^{\prime }):\,x^{\prime }\in
\Sigma ^{\prime }\}$, so that 
\begin{equation*}
T(\Sigma )=\{(y^{\prime },v(y^{\prime })):\, y^{\prime } \in \tilde{T}(\Sigma
^{\prime })\}.
\end{equation*}
Since $A_{i}^{\prime }(x_{i})\geq 1$, ($i=1,\ldots ,N-1$, $x^{\prime }\in
S^{\prime }$), we find 
\begin{eqnarray}
& & \int_{\Sigma }\varphi (x)\, {\mathcal{H}}_{N-1}(dx)  \notag
\label{basicineq1} \\
&=&\int_{\Sigma ^{\prime }}\exp \left\{ -\sum_{i=1}^{N-1}\frac{
A_{i}(x_{i})^{2}}{2}-\frac{u(x^{\prime })^{2}}{2}\right\}
\prod_{i=1}^{N-1}A_{i}^{\prime }(x_{i})\sqrt{1+\sum_{i=1}^{N-1}u_{x_{i}}(x^{
\prime })^{2}}\,dx^{\prime }  \notag \\
&=&\int_{\tilde{T}(\Sigma ^{\prime })}\exp \left\{ -\frac{\left\vert
y^{\prime }\right\vert ^{2}}{2}-\frac{v(y^{\prime })^{2}}{2}\right\} \sqrt{
1+\sum_{i=1}^{N-1}v_{y_{i}}(y^{\prime })^{2}[A_{i}^{\prime
}(A_{i}^{-1}(y_{i}))]^{2}}\,dy^{\prime }  \notag \\
&\geq &\int_{\tilde{T}(\Sigma ^{\prime })}\exp \left\{ -\frac{\left\vert
y^{\prime }\right\vert ^{2}}{2}-\frac{v(y^{\prime })^{2}}{2}\right\} \sqrt{
1+\sum_{i=1}^{N-1}v_{y_{i}}(y^{\prime })^{2}}\,dy^{\prime }  \notag \\
&=&\int_{T(\Sigma )}e^{-|x|^{2}/2}\, {\mathcal{H}}_{N-1}(dx).
\end{eqnarray}
Next assume that
$S\cap \partial M$ is a finite, disjoint union of graphs $\Sigma _{k}$ as in
(\ref{defsigma}), and of a compact set $U$ whose projection into the $x^{\prime
}$-hyperplane has ${\mathcal H}_{N-1}$-measure zero, 
\begin{equation}
S\cap \partial M=U\cup \bigcup_{k}\Sigma _{k}.  \label{union}
\end{equation}
Clearly we have 
\begin{equation}
\int_{U}\varphi (x)\,{\mathcal{H}}_{N-1}(dx)=\int_{T(U)}e^{-|x|^{2}/2}\, {
\mathcal{H}}_{N-1}(dx).  \label{basicineq2}
\end{equation}
Using (\ref{basicineq1}) and (\ref{basicineq2}) we find 
\begin{eqnarray*}
P_{\mu }(M,S) &=&\sum_{k}\int_{\Sigma _{k}}\varphi (x)\, {\mathcal{H}}
_{N-1}(dx)+\int_{U}\varphi (x)\, {\mathcal{H}}_{N-1}(dx) \\
&\geq &\sum_{k}\int_{T(\Sigma _{k})}e^{-|x|^{2}/2}\, {\mathcal{H}}
_{N-1}(dx)+\int_{T(U)}e^{-|x|^{2}/2}\, {\mathcal{H}}_{N-1}(dx) \\
&=&(2\pi )^{N/2}P_{\gamma _{N}}\left( T(M)\right) .
\end{eqnarray*}
If $M$ is a smooth subset of $S$, we can approximate it by sets satisfying (\ref{union}), 
so that inequality (\ref{perineq}) holds in the general case, too. 

Finally assume that equality holds in \eqref{Pmuisop}. Then we have 
$$
P_{\mu }(M,S)=(2\pi )^{N/2}P_{\gamma _{N}}\left( T(M)\right) 
$$
and
$$
P_{\gamma _{N}}\left( T(M)\right) =P_{\gamma _{N}}\left( H_\lambda\right) .
$$ 
Since the Gaussian isoperimetric inequality is achieved only for half-space, modulo a rotation, 
we deduce that $T(M)$ is a half-space. Hence the conclusion follows  by the definition of $T$.
$\hfill \Box $ 
\\[0.1cm]
In view of Lemma \ref{Lemma1} one has the following

\begin{corollary}\label{Corollary1}
The conclusion of Theorem \ref{Theorem2} holds for measures 
$\mu $ like (\ref{measure2}), that is
\begin{equation*}
d\mu (x)=\exp \left\{ -\frac{|x|^{2}}{2}-\sum_{i=1}^{N-1}B_{i}(x_{i})\right
\} \,dx,  
\end{equation*}
where $B_{i}\in C^{2}(a_{i},b_{i})$ with $B_{i}^{\prime \prime }(x_{i})\geq 0
$ on $(a_{i},b_{i})$, ($i=1,\ldots ,N-1$). 
\end{corollary}

{\sl Proof: } Define
$A_i \in C^3 (a_i , b_i )$ by 
$$
A_i (y):= E^{-1} \left( c_i \int _{a_i } ^y e^{-t^2 /2 - B_i (t) } \, dt \right), \ \ a_i <y<b_i ,
$$
where 
$$
c_i := \frac{1}{ \int_{a_i } ^{b_i } e^{-t^2 /2 -B_i (t)} \, dt } ,\ \ i=1, \ldots N-1.
$$
Then
$$
\prod_{i=1} ^{N-1} A_i ' (x_i ) \, 
\mbox{exp} \left\{ -\frac{x_N ^2 }{2} -\frac{1}{2} 
\sum_{i=1} ^{N-1} A_i (x_i ) ^2 \right\} = (2\pi ) ^{N-1}  \prod_{i=1} ^{N-1} c_i \, \mbox{exp} 
\left\{ - \frac{|x|^2 }{2}-\sum_{i=1} ^{N-1} B_i (x_i ) \right\} ,
$$
and by Lemma \ref{Lemma1} we have 
$$
A_i ' (y) \geq 1 \ \ \mbox{ for } \ a_i <y<b_i , \ \ i= 1, \ldots , N-1 .
$$
Now the assertion follows from Theorem \ref{Theorem2}.
$\hfill \Box $

\begin{remark} \rm  {\bf (a)} Assume that $A\in C^3 (I)$ is given and satisfies (\ref{Alimits}) and (\ref{A'geq1}), 
and define a function $B\in C^2 (a,b)$ by (\ref{formulaA}) with $c=1$. Such assumptions, as the following example shows, 
do not imply that $B''(x)\geq 0 $ on $(a,b)$.
\\
Let $A(x):= x+ \alpha x^3 $, with $\alpha >0 $, and $a=-\infty $, $b=+\infty $. Then a short computation shows that 
$$
B(x) = \alpha x^4  + \alpha ^2 x^6 /2 - \log (1+ 3 \alpha x^2 ) + \log \sqrt{2\pi } ,
$$
that is, 
$$
B''(x) = 12 \alpha x^2 + 15 \alpha ^2 x^4 + \frac{-6\alpha + 18 \alpha ^2 x^2 }{(1+ 3\alpha x^2 )^2 } .
$$
Hence $B''(0) =-6 \alpha <0$. 

\noindent   {\bf (b)} The above example also shows that Theorem \ref{Theorem2} does not follow from Corollary \ref{Corollary1}.

\noindent   {\bf (c)} In the case that $\mu $ is the Gauss
measure $\gamma _N $ and that $\omega $ is convex, it has been proved that 
$\kappa _1 \geq 1 $, see  \cite{AN} and \cite{BranChHTro}.
Together with Theorem \ref{Theorem1}, this suggests the following
conjecture:
\\
{\sl If $\omega $ is convex, then the sets $S_{\lambda }$, with $\lambda \in \mathbb{R}, $ are minimizers in the isoperimetric inequality for 
$\Omega $.}
\end{remark}

\section{Applications}

For sake of simplicity we consider the measure $\mu $ defined by 
\eqref{measure2}. We need some notation. Throughout this Section $G $
will denote a smooth domain in  $S$. We will denote by 
$C_\mu$ the constant
\begin{equation}
C_{\mu }=\int_{S^{\prime }}\exp \left( -\dfrac{\left\vert x^{\prime
}\right\vert ^{2}}{2}-\sum_{i=1}^{N-1}B_{i}(x_{i})\right) dx^{\prime }.
\label{C_zero_def}
\end{equation}
We will use the function $F:\mathbb{R}\rightarrow \mathbb{R}_{+}$ defined by 
\begin{equation}
F(t)=\int_{t}^{+\infty }\exp \left( -\frac{\sigma ^{2}}{2}\right) d\sigma ,
\label{err_funct}
\end{equation}
for any $t\in \mathbb{R}$. Such a function is strictly decreasing and
belongs to $C^{\infty }(\mathbb{R})$; we will denote by $F^{-1}:\mathbb{R}
_{+}\rightarrow \mathbb{R}$ its inverse function.

If $\Gamma $ is an open portion of 
$\partial \Omega $ with ${\mathcal H}  _{N-1} (\Gamma )>0$, let 
 $W_{\Gamma}(\Omega ,d\mu )$ be the weighted Sobolev space consisting in the 
 set of all weakly differentiable functions $u$ satisfying the following conditions: 
\begin{equation}
\left\Vert u\right\Vert _{W_{\Gamma }(\Omega ,d\mu )}^{2}:=\int_{\Omega
}\left\vert Du\right\vert ^{2}d\mu +\int_{\Omega }\left\vert u\right\vert
^{2}d\mu <+\infty ;  
\label{norma}
\end{equation}
\begin{equation*}
\left\{ 
\begin{array}{lll}
& \hbox{ there exists a sequence of functions }u_{n}\in C^{1}(\overline{
\Omega })\>\hbox{such that  } &  \\ 
& u_{n}(x)=0
\hbox{ on }\ \Gamma \
\hbox{ and} &  \\ 
& \displaystyle\lim\limits_{n\rightarrow \infty }\left( \int_{\Omega
}\left\vert D\left( u_{n}-u\right) \right\vert ^{2}d\mu +\int_{\Omega
}\left\vert u_{n}-u\right\vert ^{2}d\mu \right) =0. & 
\end{array}
\right. 
\end{equation*}
The space $W_{\Gamma }(\Omega ,d\mu )$ will be endowed with the norm defined by 
\eqref{norma}.

Now we recall a few definitions and properties about weighted rearrangements.
For exhaustive treatment on this subject we refer, e.g., to  \cite{CR}, 
\cite{Ka} and \cite{Rs}.

Let $u$ be a Lebesgue measurable function defined in $ G $. Then the
distribution function of $u$ with respect to $d\mu $ is the function $
m_{u}:[0,$\hbox{ess sup}$\left\vert u\right\vert [\rightarrow \lbrack 0,\mu
\left( G \right) [$ defined by 
\begin{equation*}
m_{\mu }(t)=\mu \left( \left\{ x\in G :\left\vert u(x)\right\vert
>t\right\} \right) ,\qquad \forall t\in \mathbb{R}_{+}.
\end{equation*}
The decreasing rearrangement with respect to $\mu $ of $u$ is the function 
$$
u^{\ast }:[0,+\infty \lbrack \rightarrow \lbrack 0, \hbox{ess sup}|u|[
$$
defined by 
\begin{equation}
u^{\ast }(s)=\inf \left\{ t\in \mathbb{R}:m_{u}(t)\leq s\right\} ,\text{ }
s\in \left] 0,\mu \left( G \right) \right] .  \label{decr}
\end{equation}
Let $G ^{\bigstar}$ be the set defined by 
\begin{equation}
\label{O*}
G ^{\bigstar }:=\{(x^{\prime },x_{N}):\,x_{N}>\lambda ,\,x^{\prime }\in
S^{\prime }\},
\quad \text{with }\lambda 
=F^{-1}\left( \frac{\mu \left(G \right) }{C_{\mu }}\right) .
\end{equation}
The weighted rearrangement of $u$ (with respect to $\mu $) is the function $
u^{\bigstar }:G ^{\bigstar }\rightarrow \lbrack 0,+\infty \lbrack $
defined by 
\begin{equation*}
u^{\bigstar }(x)=u^{\ast }\left( C_{\mu }F\left( x_{N}\right) \right)
,\qquad \forall x\in G ^{\bigstar },
\end{equation*}
where $F$ is the function given by (\ref{err_funct}) and $C_{\mu }$ is the
constant defined by (\ref{C_zero_def}). 
Observe that by definition $u^{\bigstar }$ depends just on one variable, it
is a decreasing function and moreover the functions $u$ and $u^{\bigstar }$
are equimeasurable. Therefore by Cavalieri's principle, we have 
\begin{equation}
\left\Vert u\right\Vert _{L^{p}(G ,d\mu )}=\left\Vert u^{\bigstar
}\right\Vert _{L^{p}(G ^{\bigstar },d\mu )}\,,\qquad 1\leq p\leq
+\infty .  \label{cavalieri}
\end{equation}
Let $\Gamma := \partial G \cap S$.

By a result contained in \cite{T}, we deduce that any
nonnegative function belonging to the space $W_{\Gamma }(G ,d\mu )$ satisfies
the following P\'{o}lya-Szeg\"{o} - type inequality.

\begin{theorem}
\label{Theorem3}
Let $u$ be a nonnegative function in $W_{\Gamma }(G ,d\mu )$. Then it holds 
\begin{equation}
\int_{\Omega }|Du|^{2}d\mu \geq \int_{G ^{\bigstar }}|Du^{\bigstar
}|^{2}d\mu .  \label{Polya_Szego}
\end{equation}

\end{theorem}

\bigskip As a consequence of the inequality \eqref{Polya_Szego} one deduces
that $W_{\Gamma }(G ,d\mu )$ is continuously embedded in $L^{2}(G ,d\mu
).$

\begin{corollary} 
\label{Corollary2}
For any function $u$ belonging to $W_{\Gamma }(G ,d\mu ),$
we have 
\begin{equation*}
\int_{G }\left\vert u\right\vert ^{2}d\mu \leq C\int_{G
}\left\vert Du\right\vert ^{2}d\mu ,
\end{equation*}
where $C$ is a positive constant depending only on $\mu (G )$.
\end{corollary}

{\sl Proof : } By using (\ref{Polya_Szego}), \eqref{cavalieri}
 and a result contained in \cite{Ma} (Theorem 1, p. 40),  one has that there exist a constant 
$K=K(\mu (G )) \in (0,+\infty)$ such that

\begin{equation*}
\dfrac{\displaystyle\int_{G }\left\vert Du\right\vert ^{2}d\mu }{
\displaystyle\int_{G }\left\vert u\right\vert ^{2}d\mu }\geq \dfrac{
\displaystyle\int_{G ^{\bigstar }}\left\vert Du^{\bigstar }\right\vert
^{2}d\mu }{\displaystyle\int_{G ^{\bigstar }}\left\vert u^{\bigstar
}\right\vert ^{2}d\mu }=\frac{\displaystyle\int_{\lambda }^{-\infty }\left( 
\frac{du^{\bigstar }}{dx_{N}}\right) ^{2}\exp \left( -\frac{x_{N}^{2}}{2}
\right) dx_{N}}{\displaystyle\int_{\lambda }^{+\infty }\left( u^{\bigstar
}\right) ^{2}\exp \left( -\frac{x_{N}^{2}}{2}\right) dx_{N}}\geq K,
\end{equation*}
for any $u \in W_{\Gamma }(G ,d\mu)$.
$\hfill \Box $
\\[0.1cm]
{\bf Remark 5.1.} 
We explicitly observe that by Corollary \ref{Corollary2} the norm defined by (
\ref{norma}) is equivalent to the norm 
\begin{equation}
 \label{normshort}
\left\Vert u\right\Vert _{W_{\Gamma }(G ,d\mu )}=\left( \int_{G
}\left\vert Du\right\vert ^{2}d\mu \right) ^{1/2}.
\end{equation}
Henceforth we will endow the space $W_{\Gamma }(G ,d\mu )$ with the norm (\ref{normshort}).

\noindent  Finally, we recall the classical Hardy inequality (see \cite{CR}, for instance).

\begin{proposition}
\label{hardy} Let $f$ be a function belonging to $L^{1}(G , d \mu )$ and $
E $ a measurable subset of $G $. Then the following inequality holds
true 
\begin{equation}
\int_{E}|f|d\mu \leq \int_{0}^{\mu (E)}f^{\ast }(r)dr.  \label{H-I}
\end{equation}
\end{proposition}

\bigskip 

\noindent Now we consider the following class of boundary value problems 
\begin{equation}
\left\{ 
\begin{array}{ll}
-\hbox{\rm div}\left( A(x)\nabla u\right) =\mbox{exp}\,\left\{ -\frac{|x|^{2}
}{2}-\sum_{i=1}^{N-1}B_{i}(x_{i})\right\} \,f(x) & \qquad \hbox{\rm
in } G , \\ 
&  \\ 
u=0 & \qquad \hbox{\rm on }\partial G \cap S\,
\end{array}
\right.   \label{Problem}
\end{equation}
where $G $ is an open connected subset of $S$, $(N\geq 2)$, 
$A(x)=\left( a_{ij}(x)\right) $ is a symmetric $(N\times N)-$matrix
with measurable coefficients satisfying 
\begin{equation}
\exp \left\{ -\frac{|x|^{2}}{2}-\sum_{i=1}^{N-1}B_{i}(x_{i})\right\}
\,\left\vert \zeta \right\vert ^{2}\leq a_{ij}(x)\zeta _{i}\zeta _{j}\leq
C\exp \left\{ -\frac{|x|^{2}}{2}-\sum_{i=1}^{N-1}B_{i}(x_{i})\right\}
\,\left\vert \zeta \right\vert ^{2},  \label{ell}
\end{equation}
for some $C\geq 1,$ for almost everywhere $x\in G $ and for all $\zeta
\in R^{N}$. Moreover we assume that 
\begin{equation*}
f\in L^{2}(G ,d\mu ).
\end{equation*}
\bigskip 
Let $\Gamma := G\cap S$. A solution to  problem (\ref{Problem}) is a function $u$ belonging to $
W_{\Gamma }(G ,d\mu )$ such that 
\begin{equation}
\int_{\Omega }A(x)\nabla u\nabla \psi d\mu =\int_{\Omega }f\psi d\mu ,
\label{wsol}
\end{equation}
for every $\psi \in C^{1}(\bar{G})$ such that $\psi =0$ on 
$\Gamma $.
\bigskip 
The following Theorem 5.2 gives a-priori estimates for problem (\ref{Problem}). More precisely, it states that 
every rearrangement invariant norm (with
respect to $\mu $) of the solution $u$ of (\ref{Problem}) can be estimated with the same norm of the
solution $v=v^{\bigstar }$ to the problem corresponding to the operator 

$$
L^{\bigstar }=-\hbox{\rm div}\left( \mbox{exp}\,\left\{ -\frac{|x|^{2}}{2}
-\sum_{i=1}^{N-1}B_{i}(x_{i})\right\} \,\right) \nabla v\Big)
$$
and the domain $G ^{\bigstar }$.
\\[0.1cm]
\noindent
{\bf Theorem 5.2.} {\sl 
Let $u$ be the solution to problem \eqref{Problem}.
Denote by $G ^{\bigstar }$ the subset of $\mathbb{R}^{N}$ defined in 
\eqref{O*} and by $v$ the function

\begin{equation*}
v(x)=v(x_{N})=
\int_{\lambda }^{x_{N}}\left[ \exp \left( \frac{\rho ^{2}}{2}
\right) \int_{\rho }^{\infty }\exp \left( -\frac{\xi ^{2}}{2}\right)
f^{\bigstar }(\xi )d\xi \right] d\rho 
\end{equation*}

which is the solution to the problem 
\begin{equation}
\left\{ 
\begin{array}{ll}
-\hbox{\rm div}\left( \exp \left( -\dfrac{|x|^{2}}{2}-\sum
\limits_{i=1}^{N-1}B_{i}(x_{i})\right) \nabla v\right) =\exp \left( -\dfrac{
|x|^{2}}{2}-\sum\limits_{i=1}^{N-1}B_{i}(x_{i})\right) f^{\bigstar } & \quad 
\hbox{\rm in }G ^{\bigstar }, \\ 
&  \\ 
v(\lambda )=0. & 
\end{array}
\right.   \label{Probl_Symm}
\end{equation}
Then 
\begin{equation}
u^{\bigstar }(x_{1})\leq v(x_{1})\text{ a.e. in }G ^{\bigstar },
\label{Point_Est}
\end{equation}
and 
\begin{equation}
\int_{G }\left\vert Du\right\vert ^{q}d\mu \leq \int_{G ^{\bigstar
}}\left\vert Dv\right\vert ^{q}d\mu ,\quad 0<q\leq 2.  \label{Gradient_Est}
\end{equation}
}

\medskip
\noindent We will omit the proof since it follows closely the lines, for instance, 
of  Theorem 1.1 in \cite{BCM2008} (see also \cite{AT}).
\\[0.1cm]
{\bf Remark 5.2. }

\noindent  {\bf (i)} 
The existence and the uniqueness of the solutions to problems 
(\ref{Problem}) and (\ref{Probl_Symm}), respectively, are an easy consequence of Lax-Milgram
Theorem and Corollary \ref{Corollary2}.

\noindent  {\bf (ii)}
Let us assume that the right-hand side $f$ satisfies the following
summability condition 

\begin{equation*}
\int_{\lambda }^{+ \infty}\left[ \exp \left( \frac{\rho ^{2}}{2}
\right) \int_{\rho }^{\infty }\exp \left( -\frac{\xi ^{2}}{2}\right)
f^{\bigstar }(\xi )d\xi \right] d\rho  <+\infty .
\end{equation*}
Then inequality \eqref{Point_Est} gives an estimate of the norm of $u$ in $
L^{\infty }(G ,d\mu )\equiv L^{\infty }(G )$, i.e. 
\begin{equation*}
\hbox{ess sup}|u|=u^{\star }(0)\leq \hbox{ess sup}|v|=v^{\star }(0)=
 \frac{1 }{C_{\mu }}
\int_{\lambda }^{\infty }\left[ \exp \left( \frac{\rho ^{2}}{2}
\right) 
\int_{0}^{C_{\mu }F(\rho )}
f^{^{\ast }}(\sigma )
d\sigma \right] d\rho  .
\end{equation*}

\noindent  {\bf (iii)}
Since the solution $v$ to problem \eqref{Probl_Symm}
depends just on the variable $x_{1}$, it solves the one-dimensional equation 
\begin{equation*}
-\frac{d }{d x_{N}}\left( \exp \left( -\dfrac{x_{N}^{2}}{2}
\right) \frac{d v}{d x_{N}}\right) =\exp \left( -\dfrac{
x_{N}^{2}}{2}\right) f^{\bigstar }\quad \hbox{\rm in } (\lambda, \infty),
\end{equation*}
with $v(\lambda )=0.$

\end{document}